\documentclass[12pt]{amsart}
\usepackage{amssymb}
\usepackage{l-auth}

\usepackage{graphicx}

\newtheorem{theorem}{Theorem}[section]
\newtheorem{lemma}[theorem]{Lemma}
\newtheorem{prop}[theorem]{Proposition}
\newtheorem{conj}[theorem]{Conjecture}

\newtheorem{cor}[theorem]{Corollary}

\theoremstyle{definition}

\newtheorem{quest}[theorem]{Question}

\theoremstyle{remark}
\newtheorem{remark}[theorem]{Remark}

\numberwithin{equation}{section}

\newcommand{\al}{\alpha}
\newcommand{\be}{\beta}
\newcommand{\Ga}{\Gamma}
\newcommand{\ga}{\gamma}
\newcommand{\sa}{\sigma}
\newcommand{\Sa}{\Sigma}
\newcommand{\da}{\delta}

\newcommand{\La}{\Lambda}

\newcommand{\Om}{\Omega}

\newcommand{\vp}{\varphi}

\newcommand{\scd}{{\mathcal D}}

\newcommand{\scr}{{\mathcal R}}

\newcommand{\sct}{{\mathcal T}}

\newcommand{\br}{{\mathbb R}}
\newcommand{\bz}{{\mathbb Z}}

\newcommand{\p}{\partial}
\newcommand{\wt}{\widetilde}

\newcommand{\ra}{\rightarrow}

\newcommand{\hra}{\hookrightarrow}

\newcommand{\bs}{\backslash}
\newcommand{\ov}{\overline}
\newcommand{\col}{\!:\!}

\newcommand{\Hom}{\operatorname{Hom}}
\newcommand{\hol}{\operatorname{hol}}

\newcommand{\rint}{\operatorname{int}}

\newcommand{\is}{\operatorname{Isom}}

\newcommand{\Mob}{\mathop{\rm M\ddot{o}b}\nolimits}

\newenvironment{pf}{\begin{trivlist}\item[]{\bf Proof:\ }}
{\mbox{}\hfill\rule{.08in}{.08in}\end{trivlist}}

\author{Boris N. Apanasov}
\crauthor{Boris Apanasov} 

\tit{Group actions, Teichm\"uller spaces and cobordisms}
\shorttit{Group actions} 

\begin{document}

\maketit

\address{Univ. of Oklahoma, Dept of Math., Norman, OK 73019, USA}

\email{apanasov@ou.edu}

\abstract{We discuss how the global geometry and topology of manifolds depend on different group actions of their fundamental groups, and in particular, how properties of a non-trivial compact 4-dimensional cobordism $M$ whose interior has a complete hyperbolic structure
depend on properties of the variety of discrete representations of the fundamental group of its  3-dimensional boundary $\partial M$. In addition to the standard conformal ergodic action of a uniform hyperbolic lattice on the round sphere $S^{n-1}$ and its quasiconformal deformations in $S^n$, we present several constructions of unusual actions of such lattices on everywhere wild spheres (boundaries of quasisymmetric embeddings of the closed $n$-ball into  $S^n$), on non-trivial $(n-1)$-knots in $S^{n+1}$, as well as actions defining non-trivial compact cobordisms with complete hyperbolic structures in its interiors. We show that such unusual actions always correspond to discrete representations of a given hyperbolic lattice from ``non-standard" components
of its varieties of representations (faithful or with large kernels of defining homomorphisms).
}

\notes{0}{
\subclass{57, 55, 53, 51, 32, 22, 20} 
}

\section{Geometric structures and varieties of representations}

In this survey we discuss how different hyperbolic group actions determine the global geometry and topology of corresponding quotient spaces having conformally flat
structures (locally modelled on the geometry of the sphere $S^n, n\geq 3$); in particular how a non-trivial compact 4-dimensional cobordism $M$ whose interior has a hyperbolic structure
depends on properties of the variety of discrete representations of the (hyperbolic) fundamental group of its 3-dimensional boundary $\partial M$. In addition to the standard conformal ergodic action of a uniform hyperbolic lattice on the round sphere $S^{n-1}$ and its quasiconformal deformations in $S^n$, we present several constructions of unusual actions of such lattices on everywhere wild spheres (boundaries of quasisymmetric embeddings of the closed $n$-ball into  $S^n$, with a dense subset of points with wild knotting, see Section 3), on non-trivial $(n-1)$-knots in $S^{n+1}$ (Section 2), as well as actions defining non-trivial compact cobordisms with complete hyperbolic structures in their interiors (Sections 3 and 4). We shall show that such unusual actions of a given hyperbolic lattice always correspond to its discrete representations from ``non-standard" connected components
of its varieties of representations (faithful or with large kernels of defining homomorphisms, see Sections 4 and 5).

For a given smooth $n$-manifold/orbifold $M$ with fundamental group $\pi_1(M)$, one can use its Teichm\"uller space $\sct(M)=\sct(M;X)$ to describe all possible
geometric structures on $M$ modeled
on some $(X,G)$-geometry.  In other words, we have a pair $(M,s)$
where $s$ is an $(X,G)$-structure on $M$ determined by a maximal $(X,G)$-atlas on
$M$.  Two $(X,G)$-structures $(M,s_1)$ and $(M,s_2)$ on $M$ are
equivalent if there exists an $(X,G)$-diffeomorphism $f\col(M,s_1)\ra(M,s_2)$
homotopic to the identity.  An equivalence class of $(X,G)$-structures on
$M$ is called a marked $(X,G)$-structure.  Then $\sct(M)=\sct(M;X)$ is
the space of marked $(X,G)$-structures on $M$.  The space $\sct(M;X)$ is covered by the
space $\scd(M)=\scd(M;X)$ of the development maps $d\col\wt M\ra X$ of
$(X,G)$-structures (defined as extensions of $(X,G)$-atlas charts on $M$), with a metrizable topology defined by uniform convergence of
the developments on compact sets.  So we have a topology on $\sct(M)$ induced by $\scd(M)$-topology.
In particular, if our $(X,G)$-geometry is hyperbolic the
above space $\sct(M)=\sct(M;H^n)$ becomes the Teichm\"uller space of marked
complete hyperbolic structures on $M$.  In the case of conformal geometry,
$(X,G)=(S^n,\Mob(n))$, we have the Teichm\"uller space $\sct(M;S^n)$ of marked
conformal (i.e. conformally flat) structures on $M$.

Due to the uniqueness of the development map $d\col\wt M\ra X$ (up to composition with $g\in G$), each deck transformation $T_{\al}\in G(\wt M, M)$ of the universal
covering $\wt M$ corresponds to a unique (up to conjugation by $g\in G$) element $g_{\al}\in G$ such that $g_{\al}d=dT_{\al}$. This defines the holonomy
representation $d_{\ast}\col \pi_1M\ra G$ and thus defines the holonomy map
 \begin{eqnarray}\label{hol}
\hol\col\sct(M,X)\ra\Hom(\pi_1(M),G)/G\,,
\end{eqnarray}
where the group $G$ acts on the variety of representations $\Hom(\pi_1(M),G)$ by
conjugation, and the image $\hol(M,s)$ of a structure $(M,s)\in\sct(M,X)$ is
the equivalence class $G\cdot d_\ast$ consisting of representations $gd_\ast
g^{-1}$, $g\in G$. Here $\Hom(\pi_1(M),G)/G$ has the quotient topology induced by
the algebraic convergence topology on the representation variety
$\Hom(\pi_1(M),G)$, where two representations of a group are close if
they are close on generators of the group.

The holonomy representation alone does not necessarily determine
an $(X,G)$-structure on $M$.  In particular, starting with a hyperbolic
$n$-manifold $M$, $n\geq 2$, one can obtain different conformal structures with the same holonomy (cf. \cite{M, G, A2}). One of them is
uniformizable, i.e. it is conformally equivalent to $\Om/\Ga$ where $\Om=d(\wt M)\subset S^n$ is a connected component of the discontinuity set of its holonomy group $\Gamma=d_ast(\pi_1(M))\subset\Mob(n)$. The other structure is non-uniformizable, with surjective development
map $d\col\wt M\ra S^n$, $d(\wt M)=S^n$.  However these different
structures have different developments due to the Lok-Thurston-Goldman holonomy theorem (\cite{A2}, Theorem 7.1) for closed $(X,G)$-manifolds $M$. This theorem shows that the holonomy map $\hol$ in (\ref{hol})
is an open map which lifts to a local homeomorphism
$\wt\hol$ of the space of development maps to the
variety of representations,
$
\wt\hol\col\scd(M,X)\ra\Hom(\pi_1(M),G)\,.
$

The map $\hol$ is not necessarily a local homeomorphism, especially in
neighborhoods of structures with unstable holonomy representations.
 However, in neighborhoods of structures with good holonomies $d_\ast$, the map $\hol$
itself is a local homeomorphism, see Johnson-Millson \cite{JM}.  In the case of
conformal structures, such a property follows from
the Sullivan stability theorem (Apanasov \cite{A2}, Theorem 7.2):

\begin{theorem}\label{Sul}  Let $G\subset\Mob(n)$ be a non-elementary convex-cocompact discrete
group, and $\{G_m\}$ a sequence of discrete groups $G_m\subset\Mob(n)$
isomorphic to $G$ which algebraically (on generators) converge to $G$.  Then,
for sufficiently large $m$, the groups $G_m$ are convex co-compact, and their
actions on the sphere $S^n$ are quasi-conformally conjugate to the action of
the group $G$.
\end{theorem}

In contrast to Riemann surfaces whose Teichm\"uller spaces of hyperbolic and conformal structures are isomorphic, in dimensions $n\geq 3$ one has a different situation due to the  hyperbolic rigidity of closed or finite-volume hyperbolic $n$-manifolds
$M=H^n/\Ga$, $n\geq 3$, $\pi_1(M)\cong\Ga\subset\is H^n$. In fact the Teichm\"uller space $\sct(M,H^n)$ of hyperbolic structures on
$M$ degenerates into a point while its Teichm\"uller space $\sct(M,S^n)$ of conformal structures and corresponding variety $\scr_n(\Ga)$ of
conjugacy classes of representations of the hyperbolic lattice $\pi_1(M)\cong\Ga\supset \is H^n $:
 \begin{eqnarray}\label{var}
\scr_n(\Ga)=\Hom(\Ga,\Mob(n))/\Mob(n)=\Hom(\Ga,\is H^{n+1})/\is H^{n+1}
\end{eqnarray}
is often non-degenerate. On the infinitesimal level, the non-triviality of this variety for $n=3$ (i.e. non-triviality of its tangent bundle represented by the group cohomology with coefficients in the Lie algebra of the group $\Mob(3)$) was first questioned by Borel-Wallach \cite{BW}, p.221-224. The affirmative answer to this question was given by Apanasov \cite{A1} who presented the first construction of a smooth curve in this variety defining a smooth nontrivial quasiconformal deformation of a uniform hyperbolic lattice $\Ga$ - see also Apanasov \cite{A2}, Johnson-Millson \cite{JM}, Kourouniotis \cite{K}, Lafontaine \cite{L} and Sullivan \cite{S1}. Among such conformal structures (different from the hyperbolic one) on our hyperbolic $n$-manifold $M$  we have the so called quasi-Fuchsian structures which correspond to
quasi-Fuchsian representations $\rho\in\scr_n(\Ga)$. The action of a quasi-Fuchsian group
$\rho\Ga$ on the $n$-sphere $S^n=\p H^{n+1}$ is quasi-conformally conjugate to
the action of the Fuchsian group $\rho_0\Ga\subset\Mob(n-1)\cong\is H^n$.   Due to Theorem \ref{Sul}
the subvariety $\scr_{qf}(\Ga)\subseteq \scr_n(\Ga)$ of all quasi-Fuchsian representations is an open connected component of the variety of discrete representations $\rho\in\scr_n(\Ga)$.

Equivalently, one can consider a conformal $n$-manifold (orbifold) $M$ and its Teichm\"{u}ller space $\sct'(M)\subset\sct(M,S^n)$
of equivalence classes of marked conformal structures $c$ with faithful holonomy $d_\ast\col\pi_1(M)\ra\Mob(n)$. Let also
$\sct^{\circ}(M)\subset\sct(M)$ be a subset of (classes of) conformal structures on
$M$ whose development maps are non-surjective, i.e. $\sct^{\circ}(M)$ is a subset of
almost uniformizable conformal structures (see \cite{A2}, Theorem 6.63).  On the
base of the stability Theorem \ref{Sul}, one can describe the subspaces $\sct^{\circ}(M)$ and $\sct'(M)$
as follows (see \cite{A2}, Cor.7.4).

\begin{cor}\label{v-open}
Let $M$ be a closed  manifold/orbifold with a uniformizable conformal structure.  Then the subset
\begin{eqnarray}
\sct_c(M)=\{[M,c]\in\sct'(M)\col d_\ast \pi_1(M)=G\ \text{is convex
cocompact}\,,\  |Z_G|<\infty \}
\end{eqnarray}
is open in the Teichm\"uller space $\sct'(M)$, and $\sct^{\circ}(M)$ is closed in
$\sct(M)$.
\end{cor}

Applying this result to the case of a conformal manifold with a convex cocompact holonomy group $\Ga$  and with a faithful
representation $d_\ast\col\Ga\ra\Mob(n)$, one can see that the connected components of the subvariety of $\scr_n(\Ga)$ consisting of (classes) discrete  representations can be characterized by the property that any two representations in such a component are conjugate by an equivariant quasiconformal map
$f\col S^n\ra S^n$. In particular if our conformal structure on $M$ is given by a complete hyperbolic metric (i.e. a hyperbolic lattice $\pi_1(M)=\Ga$ acts on a round $n$-ball $B^n\subset S^n$), one has the space $\sct_c(M)$
of equivalence classes of marked conformal structures on $M$ whose holonomy
representations $d_\ast\col\Ga\ra\Mob(n)$ are faithful and $d_\ast\Ga$ are
convex co-compact.  Let also $\sct_{qf}(M)$ be the subspace of $\sct_c(M)$
consisting of classes of quasi-Fuchsian structures.  Such a quasi-Fuchsian
structure on $M$ is characterized by the property that its development
$d\col\wt M=H^n\ra S^n$ is the composition of the inclusion $H^n\cong
B^n(0,1)\subset\ov{\br^n}=S^n$ and a $\Ga$-equivariant quasiconformal map
$f\col S^n\ra S^n$ so that the holonomy group $d_\ast\Ga\subset\Mob(n)$ is a
quasi-Fuchsian group, $d_\ast\Ga=f\Ga f^{-1}$.

We note that here we have an ergodic action of our hyperbolic lattice $\pi_1(M)\cong\Ga$ on topologically trivially embedded spheres $S^{n-1}\hookrightarrow S^n$ (these quasi-spheres  split $S^n$ into two quasi-balls). Later we shall consider such ergodic actions on non-trivially embedded spheres as well as actions given by non-faithful representations.

\section{Hyperbolic action on non-trivial knots}

Now, instead of Teichm\"uller spaces $\sct(M)$ of closed hyperbolic n-orbifolds $M=H^n/\Ga$ and corresponding varieties $\scr_n(\Ga)$ for their fundamental groups $\pi_1(M)\cong\Ga\subset\is H^n$ in (\ref{var}), we may consider
conformal manifolds fibered over $M$ and associated varieties of representations $\scr_m(\Ga)$, $m\geq n$,
of their holonomy groups.  It is possible to interpret these varieties of
representations as the spaces of uniformizable $m$-dimensional conformal
structures on the fundamental group $\Ga\subset\is H^n$, $n\geq 2$, of
the base $M$ of such a fibration, see Apanasov \cite{A2}, Section 7.6.

Namely, a conformal $m$-dimensional
structure on the group $\Ga\subset\is H^n$, $m\geq n$, is determined by a pair $\{N,\vp\}$
where $N$ is a conformal $m$-manifold with a non-surjective development
$d\col\wt N\ra S^m$  of its universal covering space $\wt N$ onto $d(\wt
N)=\Om\subset S^m$ and $\vp$ is a monomorphism, $\vp\col\Ga\ra\pi_1(N)$,
corresponding to the short exact sequence:
\begin{eqnarray} \label{seq}
0\ra\pi_1(\Om)\ra\pi_1(N)\ra\Ga\ra 0\,.
\end{eqnarray}

Here (due to uniformizable structures \cite{A2}) we may assume
that the development map $d$ and the natural projection
$\pi\col\Om\ra\Om/G=N$ are covering maps which factor through the universal
projection $\wt N\ra N$.  We say that two pairs $\{N_0,\vp_0\}$ and
$\{N_1,\vp_1\}$ determine the same conformal structure on $\Ga$ if there
is an orientation preserving conformal homeomorphism $f\col N_0\ra N_1$
such that $f_\ast\vp_0$ and $\vp_1$ differ (up to the isotropy subgroup
$Z(\rho_0)$ of the inclusion $\rho_0\col\Ga\subset \Mob (m)$) by an
inner automorphism of $\Ga$.  The set $\sct_m(\Ga)$ of equivalence
classes $[N,\vp]$ is called the space of conformal $m$-dimensional
structures on a given group $\Ga$, and its corresponding variety of conjugacy classes of group representations is $\scr_m(\Ga)$.

As an example, we consider the fundamental group $\Ga=\pi_1(S_g)$ of a
hyperbolic surface $S_g$ of genus $g>1$.  Then the space
$\sct_3(\Ga)$ of conformal 3-structures on $\Ga$ is determined by Seifert
fibrations $N$ over the surface $S_g$ (for fibrations with nontrivial Euler numbers,
see Gromov-Lawson-Thurston \cite{GLT} and \cite{A2}, Corollary 6.77).  Such conformal 3-manifolds $N$ have either $H^2\times\br$- or
$\wt{SL_2\br}$-geometries.
Here we only notice that the Teichm\"uller spaces $\sct_3(\Ga)$ have
many connected components. This fact is based on the existence of Seifert
fibrations with non-zero Euler classes and on the
topology of $\Ga$-actions on knots in $S^3$.

Now for a uniform hyperbolic lattice
$\Ga\subset\is H^n$, we consider the problem of connectedness of the variety
$\scr_m(\Ga)$, $m>n$, of conjugacy classes of faithful discrete
representations $\rho\col\Ga\ra\Mob(m)$.  As we observed, this variety of discrete representations
contains the space $\sct_m(\Ga)$ of $m$-dimensional uniformizable
conformal structures on $\Ga$ which satisfy the exact sequence (\ref{seq}).
Obviously, the Teichm\"uller space $\sct_n(\Ga)$ is a subspace of
$\sct_m(\Ga)$.  Nevertheless, the possible non-connectedness of
the variety $\sct_n(\Ga)$ (see the next section) does not imply
non-connectedness of the bigger variety $\sct_m(\Ga)$.  In particular the topological
obstruction for connectedness of $\sct_3(\Ga)$  (nontrivial knotting
of the limit 2-sphere $\La(G)\subset S^3$ - see the next section) is
not an obstruction in dimensions $m\geq 4$.  In fact, the nerve of that
knotting $S^2\hra S^3$ is 1-dimensional, and hence the topological 2-sphere $\La(G)\subset S^3$ nontrivially knotted in $S^3$
 is unknotted in $S^m$, $m\geq 4$.

However, we shall show that in general the varieties
$\sct_n(\Ga)$ are non-connected.  To do that we shall first
establish a link between components of these varieties and $(n-2)$-dimensional
knots in the n-sphere $S^n$ which carry conservative dynamics of
conformal actions of a given uniform hyperbolic lattice
$\Ga\subset\is H^{n-1}$, see Apanasov \cite{A2}, Theorems 7.55 and 7.60:

\begin{theorem}\label{w-knot}  For a given nontrivial ribbon $(n-2)$-knot
$K\subset S^n$, $n\geq 4$, there exists a discrete faithful
representation $\rho\col\Ga\ra \Mob(n)$ of a uniform
hyperbolic lattice $\Ga\subset\is H^{n-1}$ such that the Kleinian group
$G=\rho\Ga$ acts ergodically on the everywhere wild $(n-2)$-knot
$K_\infty=\La(G)\subset S^n$ obtained as an infinite compounding of the
knot $K$, $K_\infty=\ldots\# K\# K \# K \# \ldots$.
Moreover, the varieties $\sct_n(\Ga)$ of conformal structures on $\Ga$ and
$\scr_n(\Ga)$ of conjugacy classes
of discrete faithful representations of $\Ga$ are not connected.
\end{theorem}

\begin{figure}
\centering
\includegraphics[width=10cm]{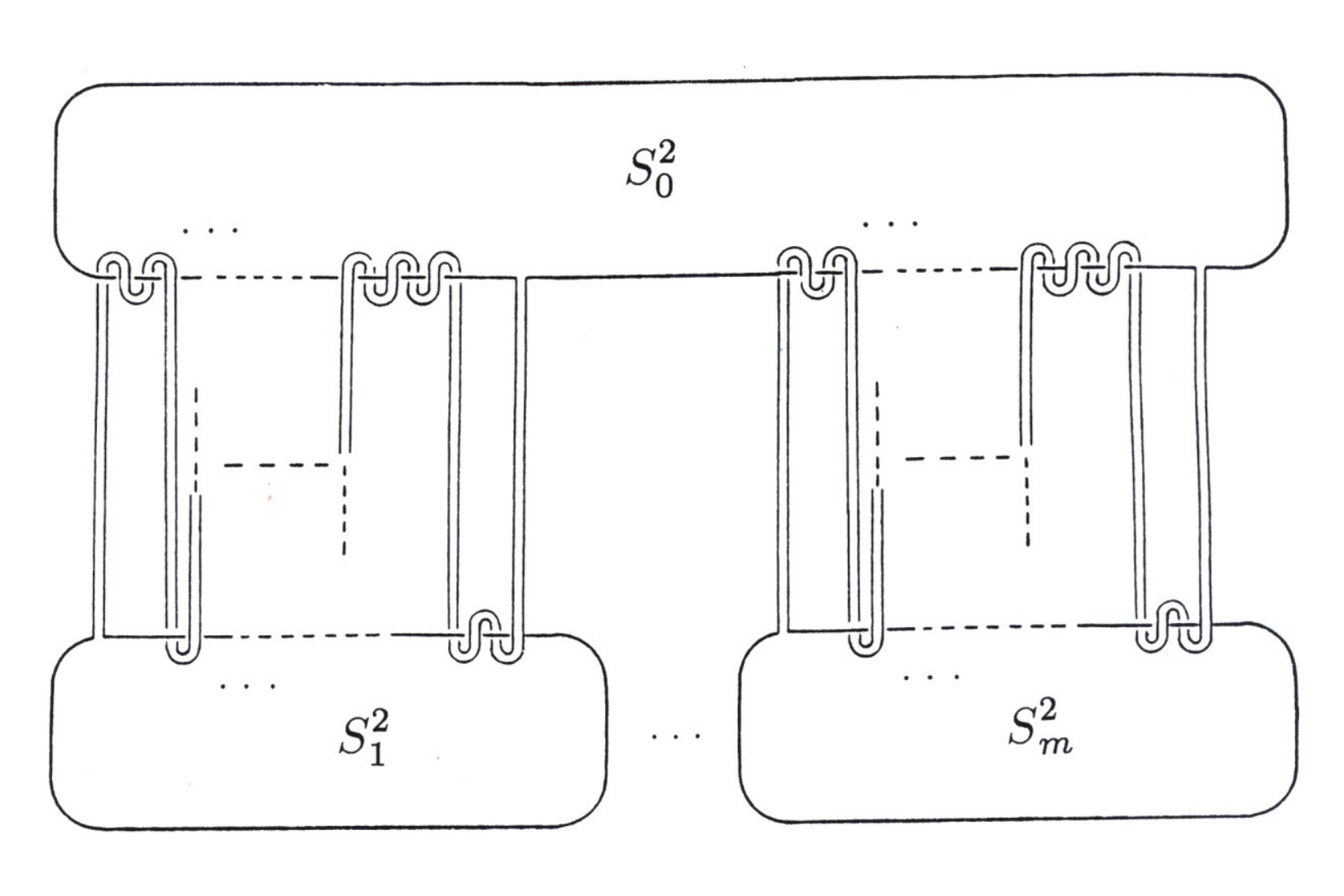}
\caption{ A ribbon 2-knot with $m$ fusions.}
\end{figure}

\begin{pf}
Here we give main steps of the proof (see details in Apanasov \cite{A2}, pp. 447-457), considering for simplicity the case $n=4$ where a nontrivial two-knot $K\subset S^4$ (not equivalent to the natural inclusion $S^2\subset S^4$) is given by an embedding $K\col S^2\hra S^4$ of
the oriented 2-sphere into the oriented 4-sphere (up to orientation preserving homeomorphisms
$f\col S^4\ra S^4$).
Simplest examples of such 2-knots can be obtained by
using the so-called suspensions and spins of classical knots in $S^3$.
The latter spun 2-knots $K\subset S^4$ can be obtained as
ribbon 2-knots. Such $(n-2)$-knots in $S^n$ generalize classical ribbon knots in $S^3$ and can be obtained as follows,
see Suzuki \cite {SS} and Figure 1.

Let $S_0\cup\ldots\cup S_m\subset\br^4$ be a trivial 2-link with $(m+1)$
components (which are trivial non-linked 2-knots) and
$f_i\col [0,1]\times B^2\hra\br^4$, $i=1,\ldots,m$, be appropriate embeddings of 3-balls, which make $m$ fusions of the 2-link.  Each of the embeddings
$f_i$ is such that
\begin{eqnarray}
f_i([0,1]\times B^2)\cap (S_0\cup\ldots\cup S_m)=f_i(\{0,1\}\times B^2)
\end{eqnarray}
has an orientation coherent with that of
the 2-link, and the disks $f_i(\{0\}\times B^2)$ and $f_i(\{1\}\times B^2)$ are
contained in different components of the link.  Then the connected sum
of the sphere $f_i(\p([0,1]\times B^2))$ and the spheres $S_0,\ldots,S_m$
represented by the homological sum
\begin{eqnarray}
(S_0\cup\ldots\cup S_m)+f_i(\p([0,1]\times B^2))=S^1_0\cup\ldots\cup S^1_{m-1}
\end{eqnarray}
is a trivial 2-link with $(m-1)$ components.
Continuing this process of fusions on the link, we finally obtain a
ribbon 2-knot with $m$ fusions, see Figure 1.

The proof of Theorem \ref{w-knot} uses the stability Theorem \ref{Sul} and is based on our block-building
method (see Apanasov \cite{A2}, Section 5.4) and geometrically controlled PL-appro\-xi\-ma\-tions of
smooth ribbon $(n-2)$-knots $K\subset S^n$ (in the conformal
category).  This  means that all spheres involved in the definition of a given ribbon knot are round (conformal) spheres in $S^n$, and each image
$f_i(B^{n-1})$ is contained in the union of finitely many round
$(n-1)$-balls $B_j$ in $S^n$, $1\leq j\leq j_i$, such that the boundary spheres of
any two adjacent balls intersect each other along a round
$(n-3)$-sphere, see Figure 2.

In other words, the $(n-1)$-dimensional ribbon $f_i(B^{n-1})$ (the union of spherical annuli) can be
obtained from a flat ribbon in $\br^{n-1}$ by sequential bendings along $(n-2)$-planes (defining bendings of corresponding discrete groups).
Here each round $(n-1)$-ball $B_j$ in the definition of our ribbon knot $K$
(either a ball from one of the ribbons $f_i(B^{n-1})$ or
one of the balls bounded by spheres $S_k$, $0\leq k\leq m$) has a discrete
action of a
hyperbolic group $G_j\subset\is H^{n-1}=\Mob(B_j)$.  Up to isotopy of the
knot $K$ and the family $\Sa$ of $(n-1)$-balls $B_j$, we may
assume that the groups $G_j$ have bending hyperbolic $(n-2)$-planes
whose boundaries at infinity $\p B_j$ are the intersection spheres $\da_j=\p
B_j\cap \p B_{j+1}$ for the adjacent balls $B_j$
and $B_{j+1}$, and that the stabilizers of $\da_j$ in $G_j$ and
$G_{j+1}$ coincide. We denote such stabilizers by $\Ga_j=G_j\cap G_{j+1}$.
This property guarantees that the amalgamated free product
\begin{eqnarray}\label{knot-g}
G=\cdots\underset{\Ga_{j-1}}\ast G_j\underset{\Ga_j}\ast
G_{j+1} \underset{\Ga_{j+1}}\ast\cdots\subset\Mob(n)
\end{eqnarray}
is a Kleinian group isomorphic to a uniform hyperbolic lattice
$\Ga\subset\is H^{n-1}$.

\begin{figure}
\centering
\includegraphics[width=10cm]{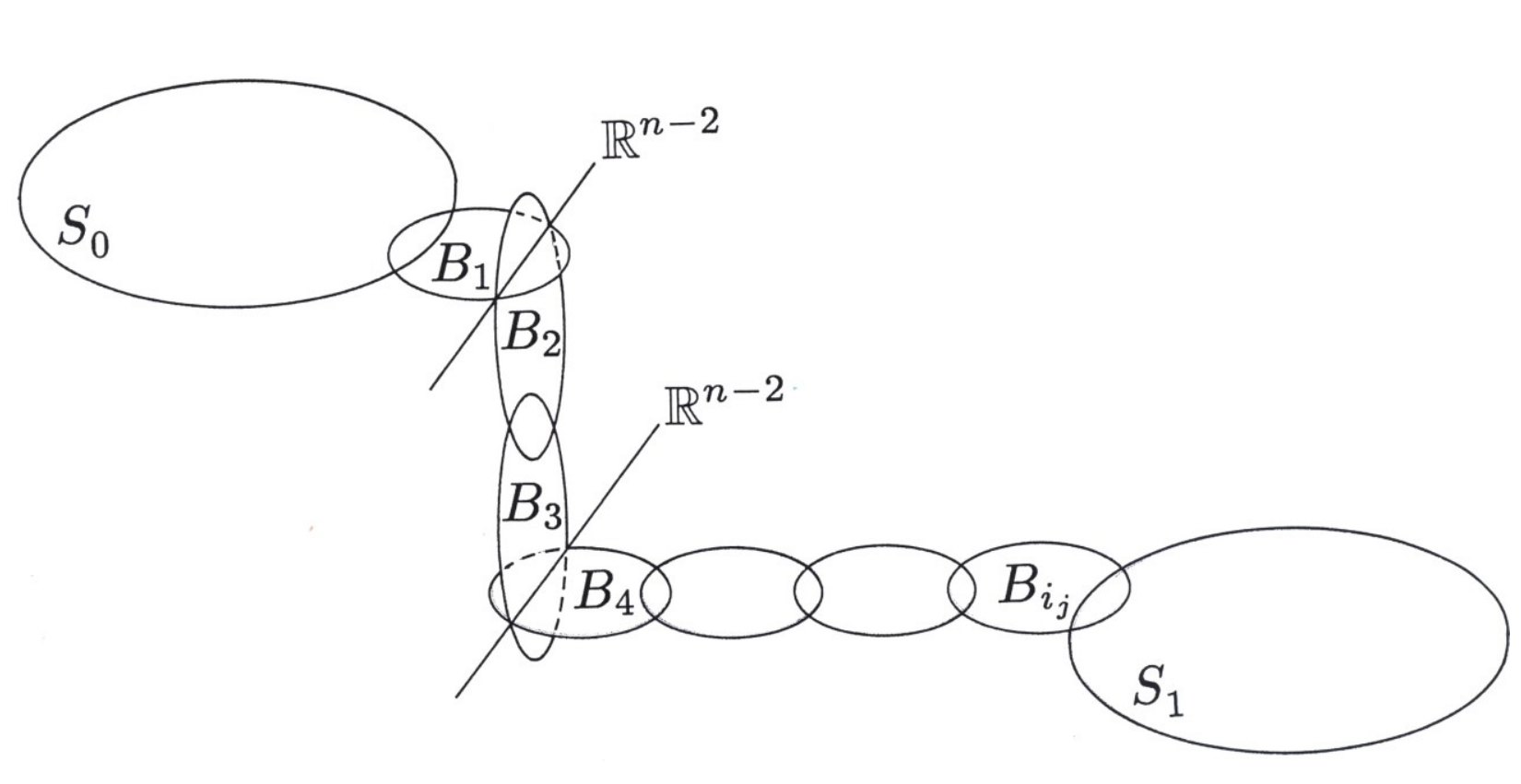}
\caption{ Block-building of a hyperbolic lattice action on $(n-2)$-knot.}
\end{figure}

This construction represents our ribbon
$(n-2)$-knot $K\subset S^n$ as the union $K_0$ of disjoint $(n-2)$-dimensional cylinders corresponding to the ribbons
$f_1,\ldots,f_m$ and $(m+1)$ disjoint round $(n-1)$-spheres with deleted disjoint round
$(n-1)$-balls (corresponding to $m$ fusions). The disjoint cylinders are the
unions of spherical $(n-2)$-dimensional annuli with disjoint interiors which lie on the boundary
spheres $\p B_j$ of the round $(n-1)$-balls $B_j$ in the construction, see Fig.2.

\begin{figure}
\centering
\includegraphics[width=10cm]{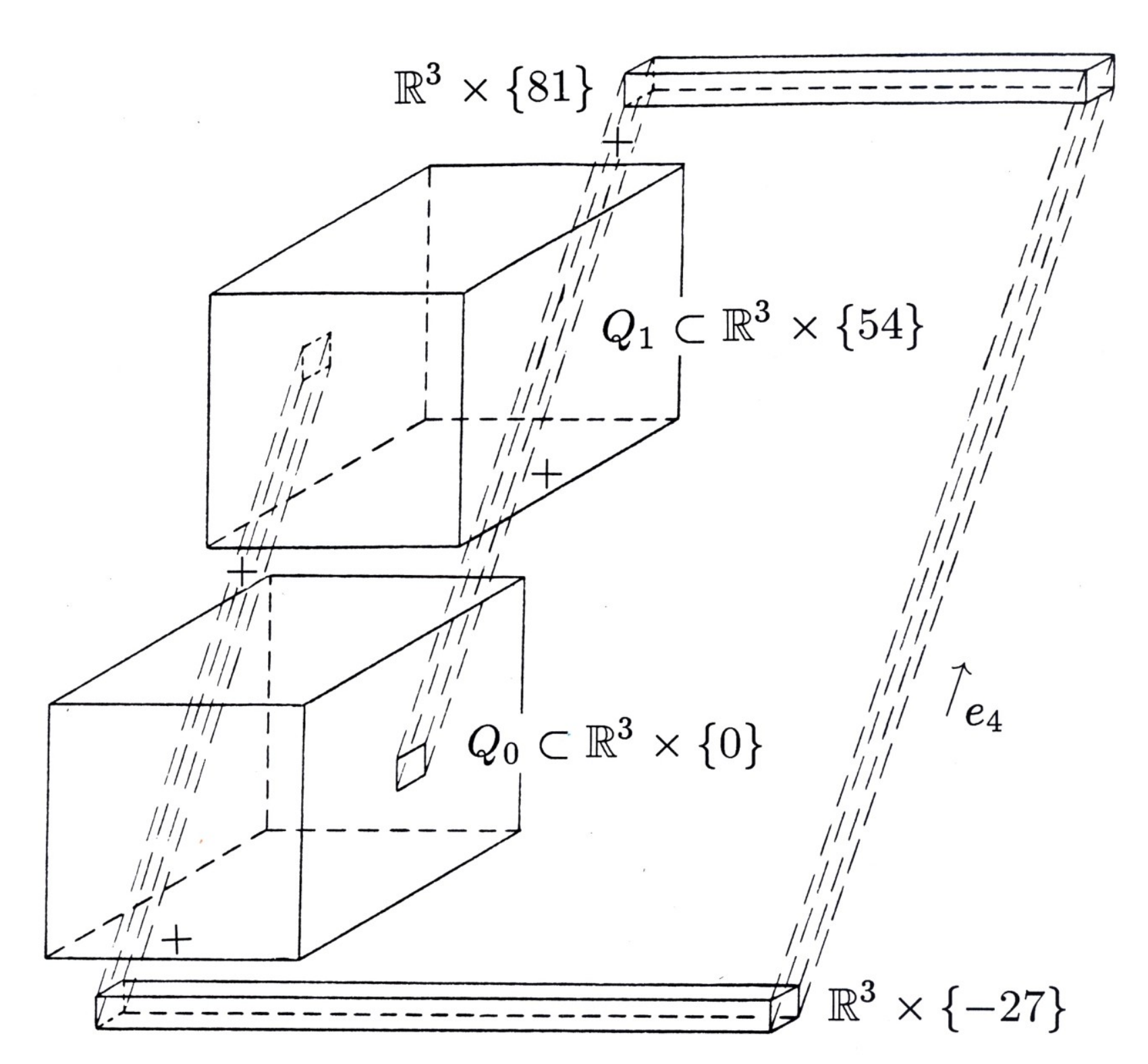}
\caption{Spun 2-knot of the trefoil as PL-ribbon knot}
\end{figure}

Using our block-groups
$G_j\subset\Mob(B_j)\subset\Mob(n)$ in (\ref{knot-g}) and sequential bendings, we obtain a uniform hyperbolic lattice
$\Ga\subset\is H^{n-1}$ having the same amalgamated free product structure as the
group $G$ in (\ref{knot-g}). The limit set $\La(G)$ of the group $G\cong\Ga\subset\is H^{n-1}$ is
the desired everywhere wild $(n-2)$-knot $K_\infty\subset S^n$ (infinite compounding of the knot
$K_0=K$), $K_\infty=\ldots\# K\# K \# K \# \ldots$. The proof of this fact is based on calculations of the Alexander invariant
of the knot $K_\infty=\La(G)$ with the complement $\Om=S^n\bs K_{\infty}$, i.e. the $\La$-module $H_\ast(\widehat\Om)$, the integral
homology $H_\ast(\widehat\Om)=H_\ast(\widehat\Om;\bz)$ with $\La$-module structure where $\La$ denotes the ring of finite Laurent polynomials
with integer coefficients, and $\widehat\Om$ is an infinite cyclic
covering space of the discontinuity set $\Om(G)$, see Apanasov \cite{A2}, Lemmas 7.57, 7.58 and Cor. 7.58.

 We illustrate how our construction
works in the simplest case of the spun 2-knot of
the classical trefoil knot $k\subset\br^3$. This 2-knot can be also represented as a PL-ribbon
knot $K\subset \br^4$ obtained by one fusion from two unlinked 2-spheres
$S_0$ and $S_1$ which are the boundaries
of 3-dimensional cubes $Q_0$ and $Q_1$ in 3-planes $\br^3\times\{0\},\,
\br^3\times\{54\}\subset \br^4$.  As the ribbon
$f_1\col B^3\hra\br^4$ we shall use the union of 3-dimensional
(smaller) cubes $Q_j$, $2\leq j\leq m$ in $\br^4$ as it is indicated in
Fig.3. Here the plus signs show the character of intersections (in our 3-dimensional projection of
$K$) of the boundaries $\p Q_0$ and $\p Q_1$ with the boundary of
3-dimensional tube which is the union $\bigcup_{2\leq j\leq m}Q_j$ of
small cubes. Edges of all cubes are parallel to the
coordinate axes in $\br^4$. Later we shall explain what is a relation between sizes of small and big cubes, as well as
 our choice of parallel
3-planes $\br^3\times\{-27\}$, $\br^3\times\{0\}$, $\br^3\times\{54\}$ and
$\br^3\times\{81\}$ in $\br^4$ which contain some of the cubes $Q_j$ and are
orthogonally joined by tubes that are boundaries of the union of the
remaining small cubes.

\begin{figure}
\centering
\includegraphics[width=9cm]{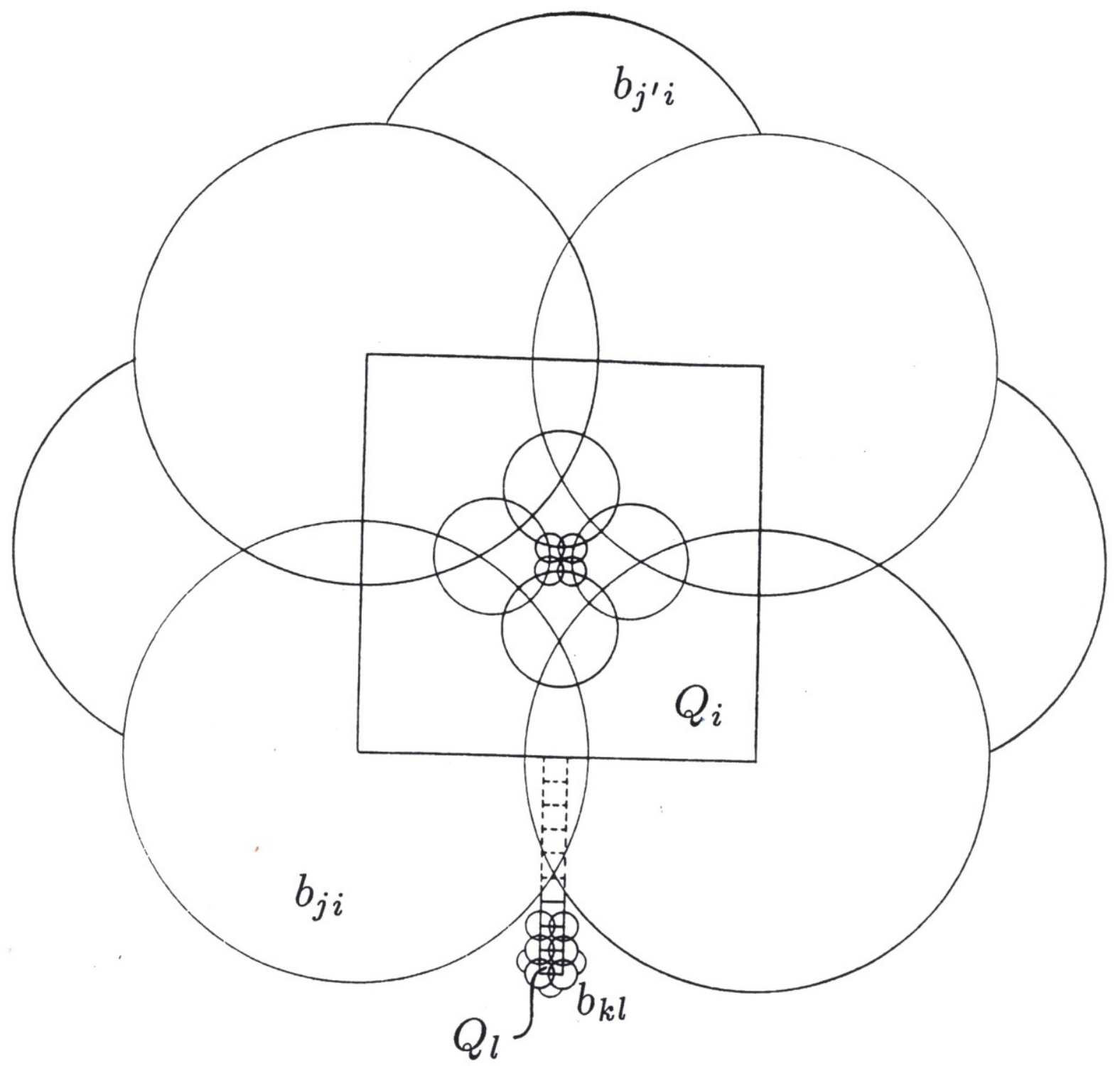}
\caption{Big and small cube sizes and ball covering}
\label{fig4}
\end{figure}

Now we define discrete block-groups $G_j$ associated with the cubes $Q_j$.
Although these $G_j$ are isomorphic to hyperbolic isometry groups in $H^3$, it is
more convenient to use quasi-Fuchsian bending deformations of these hyperbolic groups (cf. Apanasov \cite{A1, A2})
so that the obtained groups $G_j$ match the cubes $Q_j$ in the following sense.
Assuming $K=K_0$, we cover the 2-knot
$K\subset\bigcup_{0\leq j\leq m}\p Q_j$ by a family
$\Sa=\{b_{ji}\}$ of closed round 4-balls $b_{ji}$ whose boundary spheres $\p
b_{ji}$ are orthogonal to $K$.  Namely, in the first step, we take
4-balls $b_{ji}$ centered at the vertices of the cubes $Q_j$, $0\leq j\leq m$,
whose radii $r_{ji}$ are equal to each other if either $j=0,1$ or $2\leq
j\leq m$.  One more condition on these radii $r_{ji}$ is that
$b_{ji}\cap b_{kl}\neq \emptyset$ only if the centers of the different balls
$b_{ji}$ and $b_{kl}$ are the ends of a common 1-edge of one of the
cubes.  In the latter case, the magnitude of the exterior dihedral angle
bounded by the spheres $\p b_{ji}$ and $\p b_{kl}$ should equal $\pi/3$.
These 4-balls $b_{ji}$ do not cover the entire knot $K$.  On each square
2-side $X\subset\p Q_j\cap K$, we have uncovered 4-gon bounded by
circular arcs.  We cover such a 4-gon by five additional 4-balls
$b_{ji}$ centered at $X$ and whose boundary spheres $\p b_{ji}$
intersect (orthogonally) only those previously constructed balls that are
centered at the vertices of $X$.  Among these five new balls, the first
four sequentially intersect each other with external
dihedral angles $\pi/3$. The fifth ball is centered at the center of $X$
and (orthogonally) intersects
only the last new four balls, see Fig.4.  After that, we still have
uncovered those two 2-sides $X_0$ and $X_1$ of the big cubes $Q_0$ and
$Q_1$ that are (orthogonally) joined by the tube $\bigcup_{2\leq j\leq
m}Q_j$.  Here we assume that $Q_0\cap Q_2$ and $Q_1\cap Q_m$ are small
squares centered at the centers of $X_0$ and $X_1$, respectively.
Furthermore, we choose the size of the cubes $Q_j$ so that $Q_j$, $j\geq 2$,
are unit cubes and the cubes $Q_0$ and $Q_1$ have
the size which matches the covering family $\{b_{ji}\}$.  In fact, the
boundary spheres of the four additional balls centered at
$\rint (X_0)$ orthogonally intersect the corresponding spheres
centered at the four vertices of the small cube $Q_2$, see Fig.4.
This completes the construction of the family
$\Sa=\{b_{ji}\}$ of 4-balls that cover the knot $K$.  The union of these
balls, $\bigcup_{i,j}\text{int}b_{ji}=N(K)$, is a regular neighborhood $N(K)$
of the PL-ribbon knot $K$.

We can take the size of cubes $Q_j$, $j\geq 2$,
arbitrarily smaller than the size of the cubes $Q_0$ and $Q_1$. To do that,
we repeat the above process of covering the sides $X_0$ and $X_1$ by
balls $b_{ji}$ where, instead of the vertices of $X_0$ (and $X_1$),
we take the centers of the four new small balls. Then each of the annuli
in $X_0\bs Q_2$ and $X_1\bs Q_m$ will be covered by
$(4+8k)$ additional balls $b_{ji}$ (for sufficiently large integer
$k\geq 0$) instead of the above four additional balls corresponding to
$k=0$.  This allows us to take the ribbon $f_1\col B^3\hookrightarrow\br^4$ as thin
as we need. This observation makes it possible to apply
our construction of the group $G\subset\Mob(4)$ in (\ref{knot-g}) from block-groups $G_j\subset\Mob(4)$ to represent an arbitrary ribbon 2-knot
$K\subset S^4$ as the knot which lie on the boundary of the union of
3-cubes similar the above cubes $Q_j$.

Now we define a discrete block-group $G_j$  in (\ref{knot-g}) associated with a cube $Q_j$,
$0\leq j\leq m$, as the group generated by reflections with respect to
all spheres $\p b_{ji}$, that is, with respect to all spheres
$\p b_{ji}$ that intersect the cube $Q_j$.  Obviously, $G_j$ is discrete
because all spherical dihedral angles with edges $\p b_{ji}\cap \p b_{jl}$ are either
$\pi/3$ or $\pi/2$.  Furthermore,
$G_j$ preserves each of the (coordinate) 3-planes $\br^3\subset\br^4$ that contain the cube
$Q_j$.  In such a 3-plane $\br^3$, the group $G_j$ can be deformed by
bendings to a Fuchsian group $G'_j$ acting in a 3-ball $B^3\subset\br^3$, see Apanasov \cite{A2}.  This
is why we can consider the groups $G_j$ as discrete subgroups in
$\is H^3$, $G_j\cong G'_j\subset\is H^3$.
Amalgama subgroups $\Ga_j=G_j\cap G_{j+1}$ in (\ref{knot-g}) are common subgroups of groups associated with any two adjacent cubes $Q_j$ and $Q_{j+1}$.
Such a group $\Ga_j$ is generated by four reflections with respect to the spheres centered at
the vertices of the square $Q_j\cap Q_{j+1}$. Now the
Maskit combinations (see \cite{A2}, Theorem 5.17) produce a Kleinian group
$G\subset\Mob(4)$ as the free amalgamated product in (\ref{knot-g}).
We remark that for each amalgamated free product
$G_j\underset{\Ga_j}\ast G_{j+1}$, we can use a bending deformation
along the hyperbolic 2-plane $H_j$ whose boundary circle $\p H_j$ is the
limit circle of the amalgama subgroup $\Ga_j$.  As a result, we get a
new hyperbolic isometry group $G'_j\subset\is H^3$ which is isomorphic
to $G_j\underset{\Ga_j}\ast G_{j+1}$.  Applying this process $m$
times, we obtain a cocompact discrete group (a uniform hyperbolic lattice) $\Ga\subset\is H^3$
isomorphic to the group $G$.

In dimension $n=4$, there is another (non-algorithmical) way to get such a
unique hyperbolic lattice $\Ga$ by
using the Andreev \cite{An1} classification of hyperbolic compact polyhedra
in $H^3$. Namely, for the constructed group $G\subset\Mob(4)$, we can take the complement of our regular neighborhood $N(K)$
of the knot $K$ to be its fundamental polyhedron $P=P(G)\subset S^4$:
\begin{eqnarray}\label{knot-pol}
P=P(G)=\ov{\br^4}\bs N(K)\,,\quad N(K)=\bigcup_{i,j}\rint b_{ji}\,.
\end{eqnarray}
The boundary $\p P$ of the polyhedron in (\ref{knot-pol}) has
the combinatorial type of $S^2\times S^1$ where the 2-sphere $S^2$ is
decomposed into the union of spherical polygons.  In fact, $\p P$ is the
union of 3-sides each of which is the annulus on a sphere $\p b_{ji}$,
i.e. each 3-side is the product of a spherical 2-polygon $D_{ji}$ and
the circle $S^1$.  The dihedral angles between such 3-sides are determined by the
corresponding 3-dimensional dihedral angles bounded by 2-spheres  $\p
b_{ji}\cap\br^3$ in the corresponding 3-planes $\br^3\subset\br^4$, so they
are either $\pi/3$ or $\pi/2$, and the Andreev conditions apply,
see \cite{An1} and \cite{A2}, Theorem 2.41.  It follows that the combinatorial type of the
4-polyhedron $P$ determines the combinatorial type of a 3-dimensional
compact hyperbolic polyhedron $P'\subset H^3$, with the same magnitudes
of dihedral angles as those for $P$.  Thus the group $\Ga\subset\is H^3$
generated by reflections in sides of $P'$ is a uniform hyperbolic
lattice isomorphic to our group $G\subset\Mob(4)$.
\end{pf}

\section{Hyperbolic 4-cobordisms}

In this section we shall show how to change the standard conformal ergodic action of a uniform hyperbolic lattice $\Ga\subset\is H^n$ on the round sphere $S^{n-1}\subset S^n$ to its unusual actions on everywhere wild spheres (boundaries of quasisymmetric embeddings of the closed $n$-ball into  $S^n$). Such actions define non-trivial compact $(n+1)$-cobordisms with complete hyperbolic structures in their interiors. Such unusual actions always correspond to discrete representations of our hyperbolic lattice representing  ``non-standard" connected components of its variety of discrete faithful representations in $\scr_n(\Ga)$ in (\ref{var}). We shall start with hyperbolic cobordisms.

 A compact $(n+1)$-dimensional cobordism $M$ whose interior has a hyperbolic structure can be identified with a Kleinian $(n+1)$-manifold $M(G)=\left(H^{n+1}\cup\Om(G)\right)/G$,
where $G\subset\is H^{n+1}$ is a convex co-compact group of hyperbolic isometries acting in the $n$-sphere at infinity $S^n=\p H^{n+1}$ by M\"obius transformations and whose discontinuity set $\Om(G)\subset S^n$ is the union of two invariant connected components $\Omega_0$ and $\Omega_1$, $\p M=\Om(G)/G$, cf. \cite{A2}.
Since similar Kleinian 3-manifolds are always homeomorphic to surface layers $S_g\times[0,1]$ of genus $g$, it is natural to ask:  to what extent  holds the analogy with the surface layer for such $(n+1)$-cobordisms given by Kleinian $(n+1)$-manifolds $M(G)$?
\begin{figure}
\includegraphics[width=12cm]{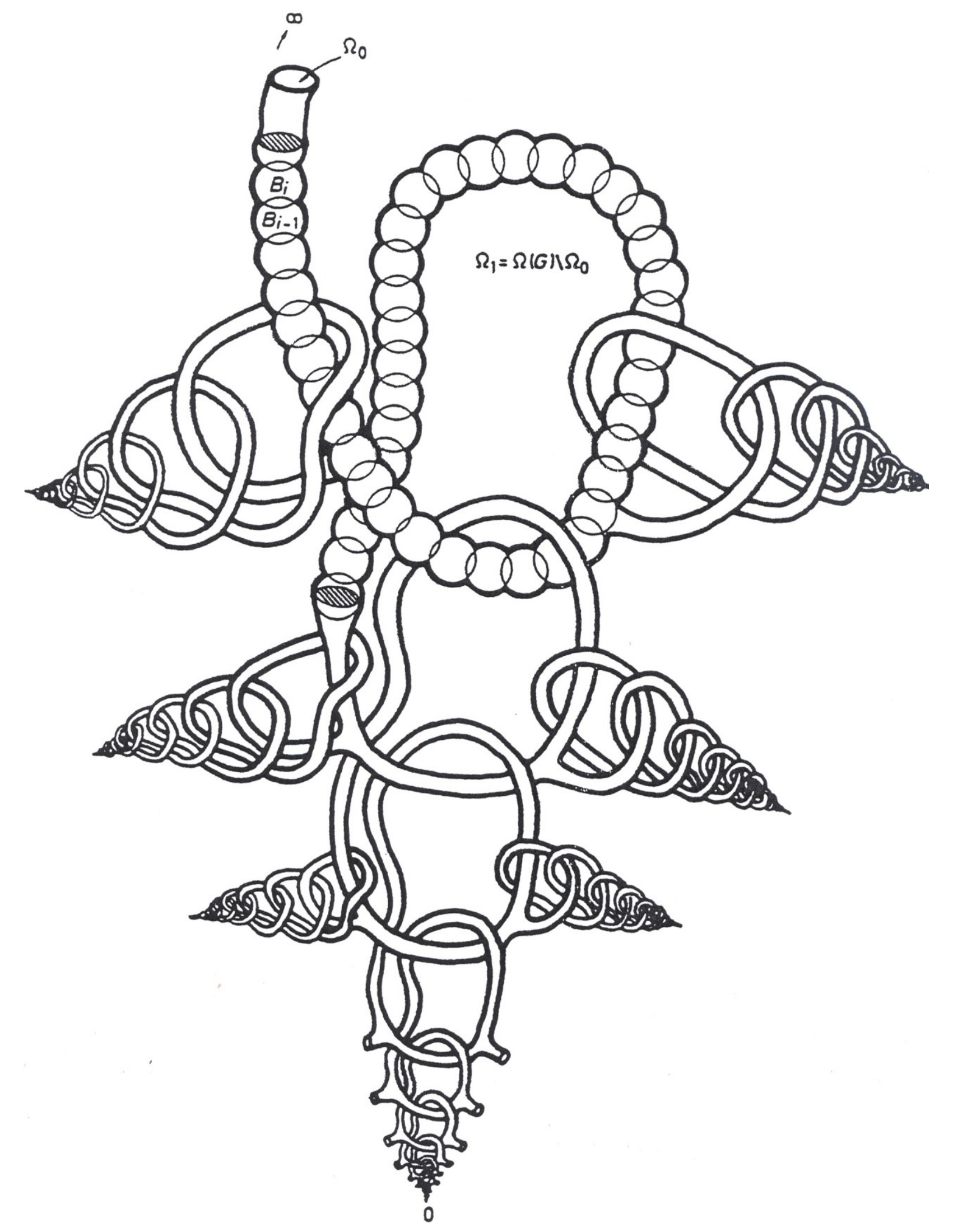}
\caption{Wild ball covering of hyperbolic 4-cobordism boundary.}
\end{figure}
One can consider such analogies of the surface layer
with various degrees of generality:
\begin{enumerate}
\item{ The product of an $n$-dimensional manifold $N_0=\Om_0/G$ and the
segment $[0,1]$.}
\item{An $h$-cobordism $(M;N_0,N_1)$ with trivial relative homotopy groups for both boundary components $N_0$ and $N_1$,
$\pi_*(M,N_0)=\pi_*(M,N_1)=0$.}
\item{A  homology cobordism $(M;N_0,N_1)$ with trivial relative  homology groups for both boundary components $N_0$ and $N_1$,
$H_*(M,N_0)=H_*(M,N_1)=0$.}
\end{enumerate}

It follows from Apanasov-Tetenov \cite{AT1} that homology cobordisms are very natural for Kleinian manifolds:

\begin{theorem} Suppose a Kleinian group $G\subset \is H^{n+1}$ is convex cocompact and has an invariant contractible component $\Om_0\subset\Om(G)\subset S^n$.
Then its compact Kleinian manifold $M(G)$
has two boundary components $N_0=\Om_0/G$ and $N_1=(\Om(G)\bs\Om_0)G$, and the triple $(M(G);N_0,N_1)$ is a homology cobordism.
\end{theorem}

It follows from the proof of this theorem (see also Cor.5.42 and Theorem 5.43 in \cite{ A2}) that an additional condition of contractibility of two invariant components
$\Om_0,\Om_1\subset \Om(G)$ of  the group $G$  turns the above homology cobordism into an
$h$-cobordism:

\begin{theorem} Let a torsion free Kleinian group $G$ acting on $S^n$, $n\geq
2$, have two invariant contractible components $\Om_0,\Om_1\subset\Om(G)$
with compact quotient manifolds $N_0=\Om_0/G$ and $N_1=\Om_1/G$.  Then the Kleinian manifold $M(G)$ is also compact, has exactly two boundary components $N_0\cup N_1=\p M(G)$, and the triple $(M(G);N_0,N_1)$ is an
$h$-cobordism,
$\pi_*(M(G),N_0)=\pi_*(M(G),N_1)=0$.
\end{theorem}

On the other hand Apanasov-Tetenov \cite{AT1} and Apanasov \cite{A4} constructions show that there are homotopically non-trivial homology cobordisms $M=M(G)$:

\begin{theorem}\label{w-ball} There exist compact Kleinian 4-manifolds $M(G)$ (both
orientable or not, and with convex cocompact hyperbolic structures in
their interiors
$\rint M(G)=H^4/G$) which are homotopically non-trivial homology cobordisms.
\end{theorem}

Convex cocompact hyperbolic isometry groups $G\subset \is H^4$ in these constructions in Theorem \ref{w-ball} have two $G$-invariant components
$\Om_0,\Om_1\subset\Om(G)\subset S^3=\p H^4$
with compact quotient manifolds $N_0=\Om_0/G$ and $N_1=\Om_1/G$. Nevertheless these boundary manifolds $N_0$ and $N_1$
differ very much from each other. While $N_0$ is homotopy equivalent to the 4-manifold $M(G)$ and is universally covered
by a quasiconformal 3-ball $\Om_0\subset S^3$, the second boundary component $N_1$ is not homotopy equivalent to $M(G)$,
and is covered by the non-simply connected $G$-invariant component $\Om_1=S^3\bs\ov{\Om_0}$ which is the complement in $S^3$ of a wildly knotted
closed 3-ball $\ov{\Om_0}\subset S^3$ shown in Fig.5. This wildly knotted closed 3-ball can be obtained by a quasisymmetric embeddings of a closed round 3-ball inextensible in neighborhoods of any boundary points (see Apanasov \cite{A3}) which can be constructed by a sequence of bendings of the initial hyperbolic 3-manifold $H^3/\Ga\approx\Om_0$ along its finitely many disjoint totally geodesic surfaces. Results of such bendings of the original hyperbolic lattice $\Ga$ are indicated in Fig.5 by sequentially intersecting balls $B_i$ forming a period of the initial Fox-Artin wild 3-ball.

So one may come to the following conjecture:

\begin{conj}\label{h-cobordism} If one had a hyperbolic 4-cobordism $M$ whose boundary components are highly (topologically and geometrically) symmetric to each other it would be in fact an h-cobordism, possibly not trivial, i.e. not homeomorphic to the product of $N_0$ and  the segment $[0,1]$.
\end{conj}

In the next section we negatively answer this Conjecture \ref{h-cobordism}. Namely we construct such hyperbolic 4-cobordisms $M=M(G)$ whose boundary components are covered by the discontinuity set $\Omega (G)\subset S^3$ with two connected components $\Omega_0$ and $\Omega_1$, where the fundamental group action $\Gamma$ is symmetric and has contractible fundamental polyhedra of the same combinatorial type allowing to realize them as a compact polyhedron in the hyperbolic 3-space, i.e.
the dihedral angle data of these polyhedra satisfy the Andreev's conditions \cite{An1}. Nevertheless we show that a geometric symmetry of boundary components of our hyperbolic 4-cobordism $M(G)$) are not enough to ensure
that the group $G=\pi_1(M)$ is quasi-fuchsian and our 4-cobordism $M$ is trivial. This is related to the
non-connectedness of the variety of discrete representations of $\Gamma$ (Teichm\"uller
space) and homomorphisms $\Gamma\rightarrow G$ with infinite kernels.

\section{Teichm\"uller Spaces and Reflection Groups in the 3-sphere}

 If the fundamental group $\Ga=\pi_1(N_0)$ of a boundary component of our hyperbolic 4-cobordism 
 $(M;N_0,N_1)$,  acts in the hyperbolic 3-space $H^3$
 as a uniform hyperbolic lattice $\Ga\subset \is H^3$, one may consider the natural conformal (conformally flat) structure on $N_0=\Om_0/G$ either as a point
 of the Teichm\"uller space $\sct(H^3/\Ga)$ of marked conformal structures on $H^3/\Ga$, or a point of the variety $\scr_3(\Ga)$ of conjugacy classes of discrete representations of the group $\Ga$ into $\is H^4$ in (\ref{var}).

Among such conformal structures on our hyperbolic $3$-manifold/orbifold $H^3/\Ga$ we have the so called quasi-Fuchsian structures which correspond to
quasi-Fuchsian representations
of the inclusion $\pi_1(N_0)\cong\Ga\subset\is H^3\subset\is H^4$ (deformations in $\scr_3(\Ga)$).  Here the Fuchsian
group $\Ga$ preserves a round ball $B^3\subset S^3=\p H^4$ and
conformally acts on this ball as a cocompact discrete group of isometries of $H^3$. Due to the Sullivan
structural stability, the
space of such quasi-Fuchsian conformal structures and the space of classes of quasi-Fuchsian representations of the hyperbolic lattice $\Ga$ into $\is H^4$
are open connected components of the Teichm\"uller space $\sct(H^3/\Ga)$ and of the variety of conjugacy classes of discrete representations
$\rho\col \Ga\ra\is H^4$,  $\scr(\Ga, \is H^4)$, respectively.
Obviously points in these (quasi-Fuchsian) components (faithful representations $\rho\col \Ga\ra \is H^4$) correspond to trivial hyperbolic 4-cobordisms
$(M(G),N_0,N_1)$ where $G=\rho(\Ga)\subset \is H^4$ has the discontinuity set $\Om(G)=\Om_0\cup\Om_1\subset S^3=\p H^4$ and $N_i=\Om_i/G$, $i=1,2$,
i.e. $M(G)$ is homeomorphic to the product of $N_0$ and the closed interval $[0,1]$.

 On the other hand, our non-trivial homology cobordisms from Theorem \ref{w-ball} correspond to non-quasi-Fuchsian discrete representations $\rho$ of $\Ga$
 which are points of another connected components consisting of discrete representations in the variety $\scr(\Ga, \is H^4)$, see Apanasov \cite{A4}.
 These discrete representations $\rho$ can be connected to the inclusion representation $i\col \Ga\subset\is H^3\subset\is H^4$ in the variety (\ref{var})
 by a continuous curve
 $\be\col [0,1]\ra \scr_3(\Ga)$, $\be (0)=i$, $\be (1)=\rho$, obtained by several deformations of the original closed hyperbolic 3-manifold $H^3/\Ga$.
 These deformations (showing non-triviality of Teichm\"uller space for a hyperbolic $n$-manifold, $n\geq 3$)  were introduced for the first time by
 Apanasov \cite{A1} and they later (after W.Thurston's Mickey Mouse example) became known as ``bendings" (of a hyperbolic 3-manifold $H^3/\Ga$ along
 disjoint totally geodesic surfaces) - see a detailed description of this curve $\be$ in Apanasov \cite{A3, A4, A2}. However it is important that
 this curve $\be([0,1]\subset\scr_3(\Ga))$ must contain some points $\be(t_0)$, $0<t_0<1$, corresponding to non-discrete representations
 $\Ga\ra \is H^4$.

Now we consider the situation of hyperbolic 4-cobordisms $M(\rho(\Ga))$ corresponding to uniform hyperbolic lattices $\Ga\subset \is H^3$ generated by reflections (or cobordisms related to their finite index subgroups). Natural inclusions of these lattices into $\is H^4$ act at infinity $\p H^4 = S^3$ as Fuchsian groups
$\Ga\subset \Mob (3)$ preserving a round ball in the 3-sphere $S^3$. In this case the  Conjecture \ref{h-cobordism} can be reformulated as the
following question on the M\"obius action of corresponding reflection groups $G=\rho(\Ga)\subset \is H^4$ on the 3-sphere $S^3=\p H^4$:

\begin{quest}\label{quest} Is any discrete M\"obius group $G$ generated by finitely many
reflections with respect to spheres $S^2\subset S^3$ and whose
fundamental polyhedron $P(G)\subset S^3$ is the union of two contractible polyhedra
$P_1, P_2\subset S^3$ of the same combinatorial type (with equal corresponding dihedral angles) quasiconformally conjugate
in the sphere $S^3$ to some Fuchsian group preserving a round ball $B^3\subset S^3$?
\end{quest}

In the next section we answer this question negatively by constructing a corresponding discrete representation
$\rho\col\Ga\ra\is H^4\cong \Mob (3)$
of a cocompact hyperbolic lattice $\Ga\subset \is H^3$ whose image $\rho(\Ga)=G\subset\is H^4\cong \Mob (3)$ will have all the properties
of the reflection group $G\subset \Mob(3)$ in Question \ref{quest}.

Necessary and sufficient conditions on a group $G$ which guarantee an affirmative answer to Question \ref{quest} can be found in \cite{AT2}.

\section{Homomorphisms of hyperbolic groups with infinite kernels}

In this section our goal is to answer Question \ref{quest} (and Conjecture \ref{h-cobordism}). In other words, one would like to decide whether
two M\"obius groups on the 3-sphere $S^3$ generated by reflections and having combinatorially similar fundamental polyhedra
(with equal corresponding dihedral angles) are
quasiconformally conjugate, that is, whether they lie in the same component of the discrete representation variety (or Teichm\"uller space of conformal structures).

Our answer to that question is negative. In fact we present a construction which implies the following (see Apanasov\cite{A6} for detailed proof):
\begin{theorem}\label{constr}
There exists a discrete M\"obius group $G\subset \Mob (3)$ on the 3-sphere $S^3$ generated by finitely many reflections such that:
\begin{enumerate}
\item Its discontinuity set $\Om(G)$ is the union
of two invariant components $\Om_0$, $\Om_1$;
\item Its fundamental polyhedron $P\subset S^3$ has two contractible components $P_i\subset\Om_i$, $i=1,2$,
having the same combinatorial type (of a compact hyperbolic polyhedron $P_0\subset H^3$);
\item For the uniform hyperbolic lattice $\Ga\subset\is H^3$ generated by reflections in sides of the hyperbolic
polyhedron $P_0\subset H^3$ and acting on the sphere $S^3=\p H^4$ as a discrete Fuchsian group  $i(\Ga)\subset \is H^4=\Mob(3)$ preserving a round ball $B^3$ (where $i\col\is H^3\subset\is H^4$ is the natural inclusion),
the group $G$ is its image under a homomorphism $\rho\col\Ga\ra G$  but it is not quasiconformally (topologically) conjugate in $S^3$ to $i(\Ga)$.
\end{enumerate}
\end{theorem}
\begin{pf} Omitting some parts of the proof (see Apanasov\cite{A6} for our detailed proof), here we concentrate on the basic for our proof
 construction of the desired M\"obius group $G\subset \Mob (3)$ generated by reflections which are defined by a finite collection $\Sa$ of reflecting 2-spheres
$S_i\subset S^3$, $1\leq i \leq N$. As the first four spheres  we consider mutually orthogonal spheres
centered at the vertices of a regular tetrahedron in
${\mathbb R}^3$. Let $B=\bigcup_{1\leq i\leq 4} B_i$ be the union of the closed balls bounded by these four spheres,
and let $\p B$ be its boundary (a topological 2-sphere) having four vertices which are the intersection points of four triples of our spheres.
Applying a  M\"obius transformation in $S^3\cong{\mathbb R}^3\cup\{\infty\}$, we may
assume that the first three spheres $S_1, S_2$ and  $S_3$
correspond to the coordinate planes $\{x\in {\mathbb R}^3\col x_i=0\}$, and
 $S_4=S^2(0,R)$ is the round sphere of some radius $R>0$ centered at the origin. The value of the radius $R$ will be determined later.

On the topological 2-sphere $\p B$ with four vertices we consider a simple closed loop $\alpha\subset \p B$
which does not contain any of our vertices and which symmetrically separates two pairs of these vertices from each other
 as the white loop does on the tennis ball shown in Figure \ref{fig6}. This loop
$\alpha$ can be considered as the boundary of a topological 2-disc $\sigma$ embedded in the complement
   $D=S^3\setminus B$ of our four balls. Our geometric construction needs a detailed
   description of such a 2-disc $\sa$ and its boundary loop $\al=\p \sa$ obtained as it is shown in
   Figure \ref{fig7}.

The desired disc $\sa\subset D=S^3\setminus B$ can be described as the boundary in the domain $D$ of the union
of a finite chain of adjacent blocks $Q_i$ (regular cubes) with disjoint interiors whose centers
lie on the coordinate planes $S_1$ and $S_2$ and whose sides are parallel to the coordinate planes.
This chain starts from the unit cube whose center lies
in the second coordinate axis, in $e_2\cdot \mathbb R_{+}\subset S_1\cap S_3$. Then our chain goes up through small adjacent cubes centered in the coordinate plane $S_1$, at some point changes its direction to the horizontal one toward the third coordinate axis, where it turns its horizontal direction by a right angle again (along the coordinate plane $S_2$), goes toward the vertical line passing through the second unit cube centered in
$e_1\cdot \mathbb R_{+}\subset S_2\cap S_3$, then goes down along that vertical line and finally ends
at that second unit cube, see Figure \ref{fig7}. We will define the size of small cubes $Q_i$ in our block chain
and the distance of the centers of two unit cubes to the origin in the next step of our construction.
\begin{figure}
\centering
\includegraphics[width=8cm]{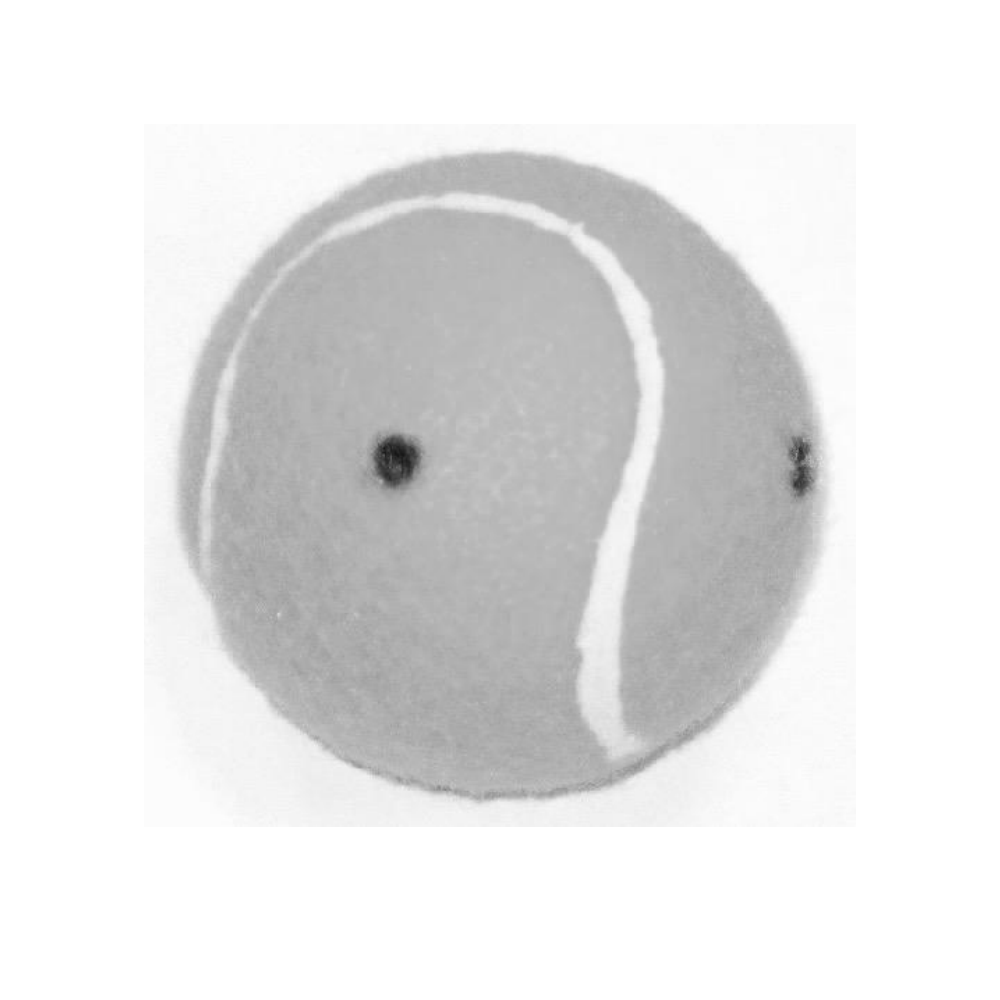}
\caption{White loop separating two pairs of vertices on a tennis ball.}
\label{fig6}
\end{figure}
\begin{figure}
\centering
\includegraphics[width=16cm]{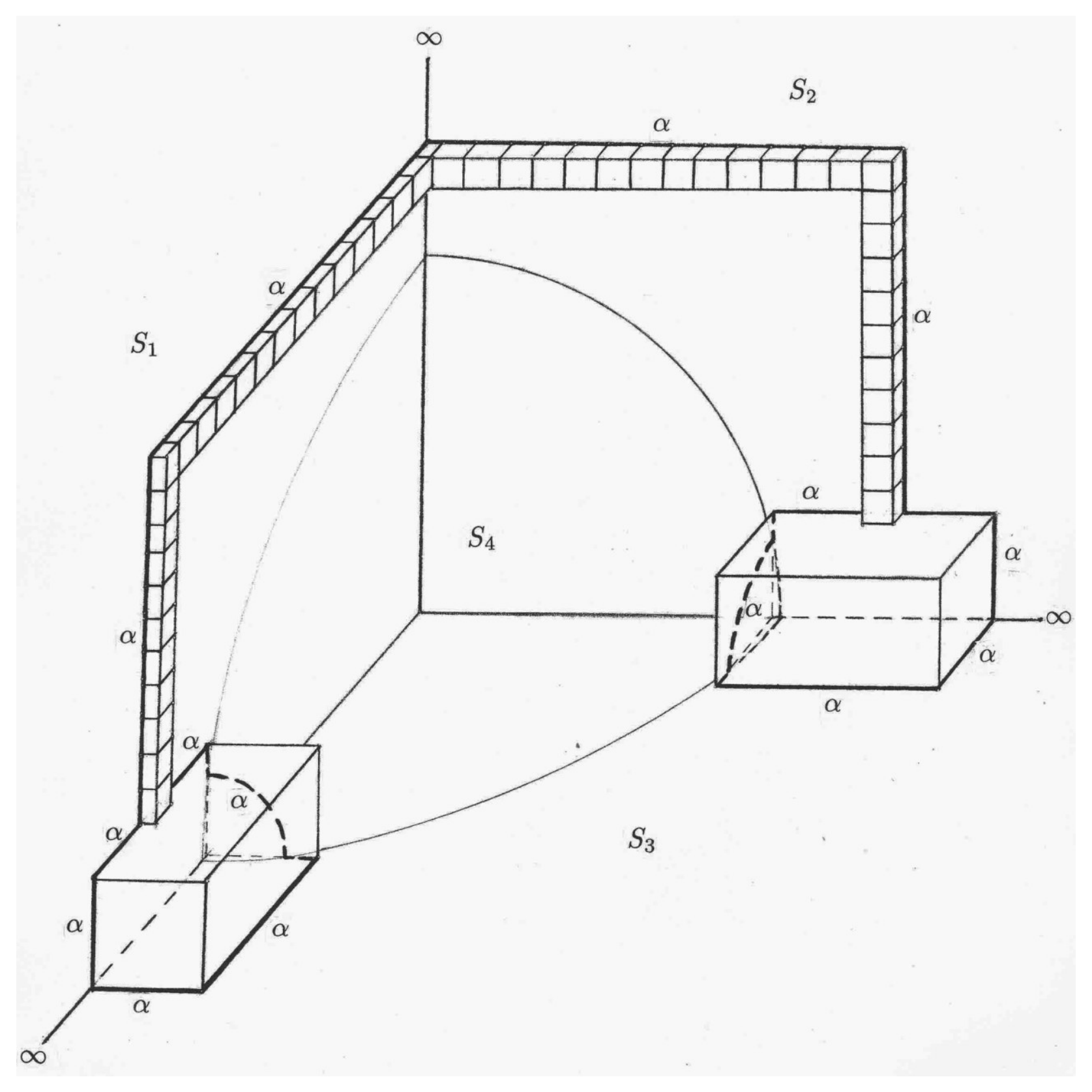}
\caption{Configuration of blocks and the loop $\alpha\subset \p B$.}
\label{fig7}
\end{figure}
Let us consider one of our cubes $Q_i$, i.e. a block of our chain, and let $f$ be its square side having a nontrivial intersection with our 2-disc
$\sa\subset D$.
For that side $f$ we consider spheres $S_j$ centered at its vertices and
having a radius such that each two spheres centered at the ends of an edge of $f$ intersect each other with angle $\pi/3$.
In particular, for the unit cubes such spheres have radius $ \sqrt{3}/3 $. From such defined spheres we select those spheres that have centers
in our domain $D$ and then include them in the collection $\Sa$ of reflecting spheres.
Now we define the distance of the centers of our big (unit) cubes to the origin. It is determined by the condition that the sphere $S_4=S^2(0,R)$
is orthogonal to the sphere $S_j\in\Sa$ centered at the vertex of such a cube closest to the origin.

As in Figure \ref{fig4}, let $f$ be a square side of one of our cubic blocks  $Q_i$ having a nontrivial intersection $f_{\sa}=f\cap\sa$ with our
2-disc $\sa\subset D$.
We consider a ring of four spheres $S_i$ whose centers are interior points of $f$ which lie outside of the four
previously defined spheres $S_j$ centered at vertices of $f$ and such that each sphere $S_i$ intersects two adjacent spheres $S_{i-1}$ and $S_{i+1}$
(we numerate spheres $S_i$ mod 4) with angle $\pi/3$. In addition these spheres $S_i$
are orthogonal to the  previously defined ring of bigger spheres $S_j$, see Figure \ref{fig4}.
From such defined spheres $S_i$ we select those spheres that have nontrivial
intersections with our domain $D$ outside the previously defined spheres $S_j$, and then include them in the collection $\Sa$ of reflecting spheres.
If our side $f$ is not the top side of one of the two unit cubes we add another sphere $S_k\in\Sa$.
It is centered at the center of this side $f$ and is orthogonal to the four previously defined spheres $S_i$ with centers in $f$, see Figure \ref{fig4}.

Now let $f$ be the top side of one of the two unit cubes of our chain. Then, as before,
we consider another ring of four spheres $S_k$. Their centers are interior points of $f$, lie outside of the four
previously defined spheres $S_i$ closer to the center of $f$ and such that each sphere $S_k$ intersects two adjacent spheres $S_{k-1}$ and $S_{k+1}$
(we numerate spheres $S_k$ mod 4) with angle $\pi/3$. In addition these new four spheres $S_k$
are orthogonal to the previously defined ring of bigger spheres $S_i$, see Figure \ref{fig4}.
 We note that the centers of these four new spheres $S_k$ are vertices of a small square
$f_s\subset f$ whose edges are parallel to the edges of $f$, see Figure \ref{fig4}. We set this square $f_s$ as the bottom side of the small cubic box
adjacent to the unit one.
This finishes our definition of the family of twelve round spheres whose interiors cover the square ring $f\bs f_s$ on the top side of one of the two unit cubes
in our cube chain and tells us which two spheres among the four new defined spheres $S_k$ were already included
in the collection $\Sa$ of reflecting spheres
(as the spheres $S_j\in\Sa$ associated to small cubes in the first step).

This also defines the size of small cubes in our block chain. Now we can vary the remaining free parameter $R$
(which is the radius of the sphere $S_4\in\Sa$) in order to make two horizontal rows of small blocks with centers in $S_1$ and $S_2$, correspondingly,
to share a common cubic block centered at a point in $e_3\cdot \mathbb R_{+}\subset S_1\cap S_2$, see Figure \ref{fig7}.

The constructed collection $\Sa$ of reflecting spheres $S_j$ bounding round balls $B_j$, $1\leq j\leq N$, has the following properties:

\begin{enumerate}
\item The closure of our 2-disc $\sa\subset D$ is covered by balls $B_j$: $\bar{\sigma}\subset \rint\bigcup\limits_{j\geq5}^{N}B_j$;
\item Any two spheres $S_j, S_{j'}\in\Sa$ 
either are disjoint or intersect with angle $\pi/2$ or  $\pi/3$;
\item The complement of all balls, $S^3\setminus\bigcup\limits_{j=1}^{N}B_j$ is the union of two contractible
polyhedra $P_1$ and $P_2$ of the same combinatorial type.
\end{enumerate}

Therefore we can use the constructed collection $\Sa$ of reflecting spheres $S_i$ to define a discrete group $G=G_{\Sa}\subset \Mob(3)$
generated by $N$ reflections in spheres $S_j\in \Sa$. The fundamental polyhedron $P=P_1\cup P_2\subset S^3$ for the action of this discrete reflection group
$G$ on the sphere $S^3$ is the union of two connected polyhedra
$P_1$ and $P_2$ which are disjoint topological balls. So the discontinuity set $\Omega(G)\subset S^3$ of $G$ consists
of two invariant connected components $\Om_0$ and $\Om_1$:

\begin{equation}\label{comp}
\Omega(G)=\bigcup_{g\in G} g(\bar{P})=\Om_0\cup\Om_1\,,\quad \Om_i=\bigcup_{g\in G} g(\bar{P_i})\,, \quad i=1,2.
\end{equation}
In fact, despite the contractibility of polyhedra $P_1$ and $P_2$ both components $\Om_0$ and $\Om_1$ are not simply connected and even are mutually linked:
\begin{lemma}\label{h-body}
The splitting of the discontinuity set $\Om\subset S^3$ of our discrete reflection group $G=G_{\Sa}\subset \Mob(3)$ into $G$-invariant components\, $\Om_0$ and $\Om_1$ in \ref{comp} defines a Heegaard splitting of the 3-sphere $S^3$ of infinite genus with ergodic word hyperbolic group $G$ action on the separating boundary $\La(G)$ which is quasi-self-similar in the sense of Sullivan.
\end{lemma}

\begin{figure}
\centering
\includegraphics[width=16cm]{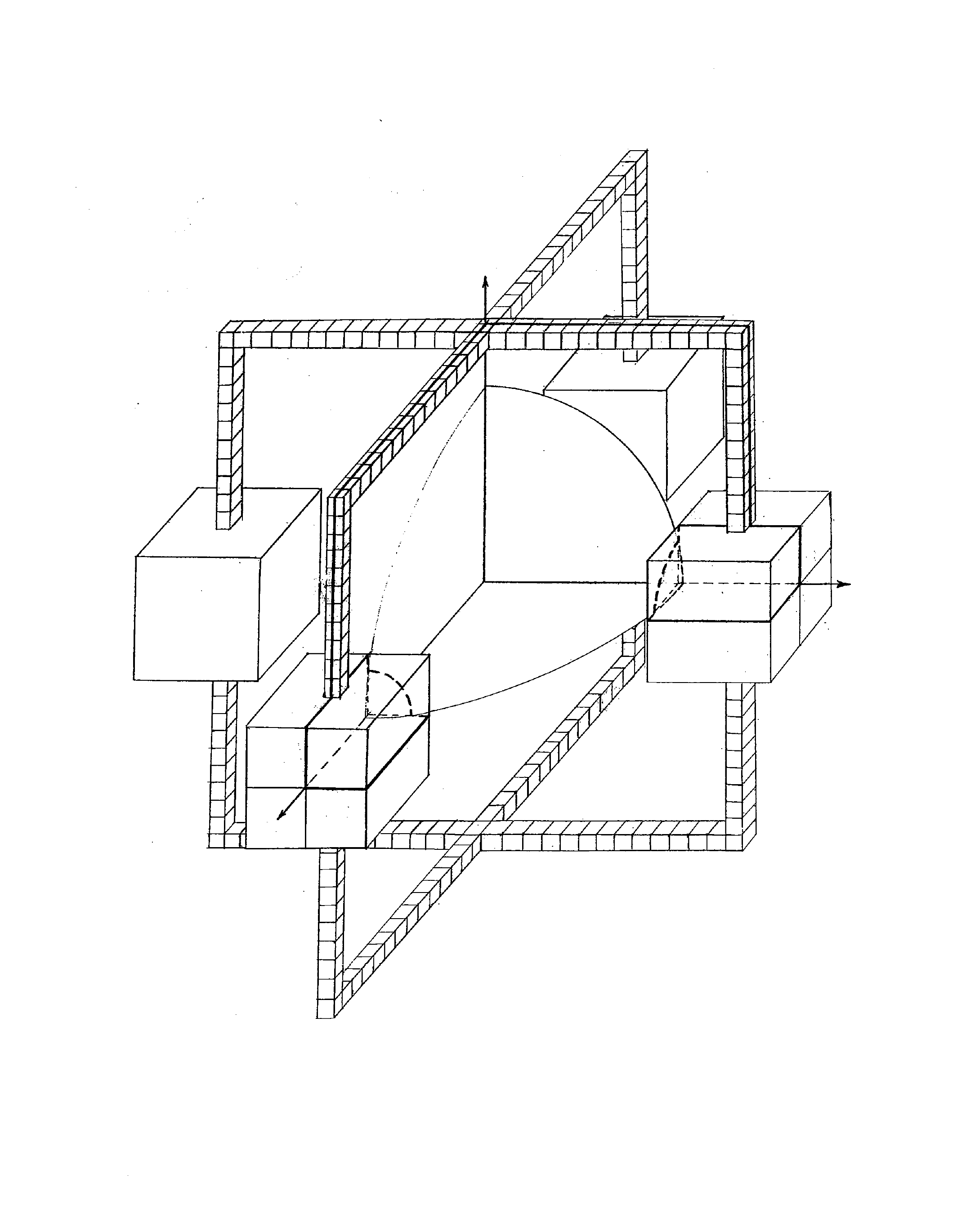}
\caption{Handlebody obtained by the first 3 reflections of the cub chain.}
\label{fig8}
\end{figure}

This fact is illustrated by Figure \ref{fig8} where one can
see a handlebody obtained from our initial chain of building blocks in Figure \ref{fig7} by the union of the images of this block chain by first generating reflections in the group $G$ (in $S_1, S_2$ and $S_3$). Then one has a non-contractible loop $\be_0\subset\Om_0$ (represents a non-trivial element of the fundamental group $\pi_1(\Om_0)$) which
lies inside of this handlebody in Figure \ref{fig8}  and is linked with the second loop $\be_1\subset\Om_1$ in the unbounded component $\Om_1$ which goes around $\bar{P_1}\cup g_3(\bar{P_1})$, i.e. around one of the handles of the handlebody in Figure \ref{fig8}. This loop $\be_1$
represents a non-trivial element of the fundamental group $\pi_1(\Om_1)$.
The resulting
handlebodies $\Omega_0$ and  $\Omega_1$ are the unions of the corresponding images $g(\bar{P_i})$ of the polyhedra $\bar{P_0}$ and $\bar{P_1}$, so they have infinitely many mutually linked handles. Their fundamental groups $\pi_1(\Om_0)$ and  $\pi_1(\Om_1)$ have infinitely many generators, and some of those generators correspond to the group $G$-images of the linked loops
$\be_0\subset\Om_0$ and $\be_1\subset\Om_1$. The limit set
$\La(G)$ is the common boundary of $\Om_0$ and $\Om_1$. Since the group $G\subset\Mob(3)$ acts on the hyperbolic 4-space $H^4$, $\p H^4=S^3$,
as a convex cocompact isometry group, its action on the limit set $\La(G)$ is ergodic. Moreover, the common boundary
$\La(G)$ of the handlebodies $\Omega_0$ and  $\Omega_1$ is quasi-self-similar in the sense of Sullivan \cite{S3}
(see Apanasov\cite{A6} for details).

Another important observation in our proof of Theorem \ref{constr} is that the combinatorial type (with magnitudes of dihedral angles) of the bounded component $P_1$ of the fundamental polyhedron $P\subset S^3$ coincides with
the combinatorial type of its unbounded component $P_2$. Applying Andreev's theorem on 3-dimensional hyperbolic polyhedra \cite{An1},
one can see that there exists a compact hyperbolic polyhedron $P_0\subset H^3$ of the same combinatorial type with the same dihedral angles ($\pi/2$ or $\pi/2$).
So one can consider a uniform hyperbolic lattice $\Ga\subset\is H^3$ generated by reflections in sides of the hyperbolic
polyhedron $P_0$. This 3-hyperbolic lattice $\Ga$ acts in the sphere $S^3$ as a discrete co-compact Fuchsian group $i(\Ga)\subset \is H^4=\Mob(3)$
(i.e. as the group
$i(\Ga)\subset \is H^4$ where $i\col\is H^3\subset\is H^4$ is the natural inclusion)
preserving a round ball $B^3$ and having its boundary sphere $S^2=\p B^3$ as the limit set. Obviously there is no self-homeomorphism of the sphere $S^3$
conjugating the action of the groups $G$ and $i(\Ga)$ because the limit set $\La(G)$ is not a topological sphere. So the constructed group $G$ is not a
quasi-Fuchsian group.

One can construct a natural homomorphism $\rho\col \Ga\ra G$, $\rho\in\scr_3(\ga)$, between these two Gromov hyperbolic groups $G\subset \is H^4$ and $\Ga\subset \is H^3$
defined by the correspondence between sides of
the hyperbolic polyhedron $P_0\subset H^3$ and reflecting spheres $S_i$ in the collection $\Sa$ bounding the fundamental polyhedra $P_1$ and $P_2$.

\end{pf}

The homomorphism $\rho$ cannot be an isomorphism since its kernel $\rho^{-1}(e_G)$ is not trivial: 

\begin{prop}\label{homo}
The homomorphism $\rho\in\scr_3(\Ga)$, $\rho\col \Ga\ra G$, in Theorem \ref{constr} is not an isomorphism. Its kernel $\ker(\rho)=\rho^{-1}(e_G)$ is a
free rank 3 subgroup $F_3\lhd\Ga$.
\end{prop}

In fact the kernel $\rho^{-1}(e_G)$ is a
free rank 3 group $F_3=\langle x, y, z\rangle$ generated by three hyperbolic translations $x, y, z \in\Ga$. The first hyperbolic translation  $x=a_1b_1$ in $H^3$ is the composition of reflections $a_1$ and $b_1$ in two disjoint hyperbolic planes $H_1, H'_1\subset H^3$ containing those two 2-dimensional faces of
the hyperbolic polyhedron $P_0$ that correspond to two sides of the polyhedron $P_1$ which are disjoint parts of the sphere $S_4$.  The second
 hyperbolic translation  $y=a_2b_2$ in $H^3$ is the composition of reflections $a_2$ and $b_2$ in two disjoint hyperbolic planes $H_2, H'_2\subset H^3$  containing those two 2-dimensional faces of the hyperbolic polyhedron $P_0$ that correspond to two sides of the polyhedron $P_1$ which are disjoint parts of the sphere $S_3$.
 And the third generator $z$ is a hyperbolic translation in $H^3$ which is $a_1$-conjugate of $y$, $z=a_1ya_1$. The fact that these hyperbolic 2-planes $H_1$ and $H'_1$ (correspondingly, the 2-planes $H_2$ and $H'_2$) are disjoint follows from Andreev's result \cite{An2} on sharp angled hyperbolic polyhedra. Restricting our homomorphism $\rho$ to the subgroup of $\Ga$ generated by reflections $a_1, a_2, b_1, b_2\in \Ga$, we can formulate its properties as the following statement in combinatorial group theory (see Apanasov\cite{A6} for its detailed proof):

 \begin{lemma}\label{ker}
Let $A=\langle a_1, a_2 \mid a_1^2,\, a_2^2,\, (a_1a_2)^2\rangle $ $\cong$ $B=\langle b_1, b_2 \mid b_1^2,\, b_2^2,\, (b_1b_2)^2\rangle $ $\cong$
$C=\langle c_1, c_2 \mid c_1^2,\, c_2^2,\, (c_1c_2)^2\rangle \cong \mathbb{Z}_2 \times \mathbb{Z}_2$,
and let $\varphi\col A\ast B\ra C$ be a homomorphism of the free product $A\ast B$ into $C$ such that $\varphi(a_1)=\varphi(b_1)=c_1$ and
 $\varphi(a_2)=\varphi(b_2)=c_2$. Then the kernel $\ker(\varphi)=\varphi^{-1}(e_C)$ of $\varphi$ is a free rank 3 subgroup $F_3\lhd A\ast B$  generated by elements
 $x=a_1b_1$, $y=a_2b_2$ and $z=a_1a_2b_2a_1=a_1ya_1$.
\end{lemma}

Therefore the configuration of reflecting spheres $S_j\subset \Sa$ shows that one can deform our discrete co-compact Fuchsian group
$i(\Ga)\subset \is H^4=\Mob(3)$ preserving a round ball $B^3\subset S^3$ into the group $G\subset \is H^4$ by continuously moving two pairs
of reflecting 2-spheres of the Fuchsian group $i(\Ga)$ corresponding to the pairs of hyperbolic planes $H_1, H'_1\subset H^3$ and $H_2, H'_2\subset H^3$
into the reflecting spheres $S_4$ and $S_3$ while keeping all dihedral angles unchanged.

Also our construction shows that the subvariety $\scr_{dis} (\Ga,\is H^4)\subset \scr_3(\Ga)$ consisting of conjugacy classes of discrete representations of a hyperbolic lattice $\Ga\subset \is H^{3}$ generated by reflections may have several connected components.
Indeed as we showed our representation  $\rho\col\Ga\ra\is H^4$ from Proposition \ref{homo} whose image is our convex cocompact discrete reflection group
$G=\rho(\Ga)\subset\is H^4$ in Theorem \ref{constr} is not quasi-Fuchsian. The stability Theorem \ref{Sul} implies that the conjugacy class defined by this representation $\rho$ does not belong to the connected (quasi-Fuchsian) component of our subvariety
$\scr_{dis} (\Ga,\is H^4)$  consisting of conjugacy classes of all quasi-Fuchsian representations
quasiconformally conjugate to the natural inclusion $i\col\Ga\subset \is H^3\subset\is H^4$.

Moreover, counting the reflecting
spheres $S_j$ of the generators of our group $G=\rho(\Ga)\subset\is H^4$ and the generators of the corresponding hyperbolic lattice $\Ga\subset \is H^{3}$ and using Theorem 7.36 in \cite{A2} (for the co-finite case, see Huling \cite{H}), one can see that dimensions of these two connected components of $\scr_{dis} (\Ga,\is H^4)$ one of which contains the representation $\rho$ and the other one contains the inclusion $i$ are different:

\begin{prop} Let $\Ga\subset \is H^{3}$ and $G=\rho(\Ga)\subset\is H^4$ be a uniform hyperbolic lattice and its reflection group image in Theorem \ref{constr}, and let $\scr_G$ and $\scr_{qf}$ be two connected components of the subvariety
$\scr_{dis} (\Ga,\is H^{4})\subset \scr (\Ga,\is H^{4})$ consisting of conjugacy classes of discrete representations quasiconformally conjugate
in $S^3$ to, correspondingly,
the representation $\rho\col\Ga\ra G$ and the inclusion $i\col\Ga\subset \is H^3\subset\is H^4$. Then, near the inclusion $i$, the variety $\scr_i$
is a smooth manifold of dimension 775 while the variety $\scr_G$ has dimension 773.
\end{prop}
\begin{remark} Our Theorem \ref{constr} shows that a hyperbolic 4-cobordism $M$ whose boundary components are highly (topologically and geometrically) symmetric to each other is not necessarily a trivial 4-cobordism, i.e. not homeomorphic to the product. This fact is related to the Novikov conjecture \cite{N} on the homotopy invariance of characteristic numbers of a manifold derived from the fundamental group. Also, since hyperbolic manifolds are of $K(\pi_1,1)$ type (their universal covers are simply connected), it indicates that such hyperbolic 4-cobordisms (and hyperbolic groups) could be among classes for which the conjecture has been verified, see \cite{FRR, G, KL}.
\end{remark}

\begin{remark} We refer to \cite{A3, A7} for applications of our  constructions of non-trivial 4-dimensional cobordisms $M$ to geometric function theory, in particular to quasiconformal, quasisymmetric and quasiregular mappings. This is related to the M.A.Lavrentiev problem, the Zorich map with an essential singularity at infinity and a quasiregular analogue of domains of holomorphy in complex analysis.
\end{remark}


\end{document}